\documentclass[12pt]{article}

\usepackage{auto-pst-pdf} 
\usepackage{pstricks,pst-node}
\usepackage{amsmath,amsthm}
\usepackage{amssymb}
\usepackage{color}
\usepackage{xspace}
\usepackage[colorlinks=true,
linkcolor=green,
filecolor=brown,
citecolor=green]{hyperref}
\def\blue{\textcolor{blue} }

\def\mag{\textcolor{magenta} }

\def\a{\ensuremath{\mathcal A}\xspace}
\def\b{\ensuremath{\mathcal B}\xspace}

\def\gl{ground level\xspace}
\def\gf{generating function\xspace}
\def\gfs{generating functions\xspace}

\def\mbf#1{\mathchoice{\hbox{\boldmath $\displaystyle #1$}}
        {\hbox{\boldmath $\textstyle #1$}}
        {\hbox{\boldmath $\scriptstyle #1$}}
        {\hbox{\boldmath $\scriptscriptstyle #1$}}} 

\begin{document}
\newtheorem{theorem}{Theorem}
\newtheorem{defn}[theorem]{Definition}
\newtheorem{lemma}[theorem]{Lemma}
\newtheorem{prop}[theorem]{Proposition}
\newtheorem{cor}[theorem]{Corollary}
\begin{center}
{\Large
\blue{A bijection for two sequences in OEIS}   \\ 
}

\vspace{2mm}

David Callan  \\

\vspace{1mm}

\today
\end{center}
\vspace*{2mm}

\begin{abstract}
We show that sequences A026737 and A111279 in The On-Line Encyclopedia of Integer Sequences are the same by giving a bijection between two classes of Grand Schr\"{o}der paths.
\end{abstract}

\section{Introduction} \label{intro}

\vspace*{-2mm}

In a comment on sequence \htmladdnormallink{A026737}{http://oeis.org/A026737} in OEIS \cite{oeis}, 
Andrew Plewe asks if it is the same as  \htmladdnormallink{A111279}{http://oeis.org/A111279}. 
The answer is yes.
As we will see, each of the sequences counts a class of \emph{Grand Schr\"{o}der paths}, 
that is, lattice paths of upsteps $U=(1,1)$, flatsteps $F=(2,0)$, and downsteps $D=(1,-1)$ 
starting at the origin (0,0) and ending on the $x$-axis, with size 
measured by \# upsteps + \# flatsteps. Thus a Grand Schr\"{o}der path of size $n$ ends at $(2n,0)$. 

The $(n+1)$-th term of A026737 counts Grand Schr\"{o}der paths of size $n$ \emph{all of 
whose flatsteps $($if any$\hspace{.1em})$ lie on the horizontal line $y=2$}. This 
assertion follows immediately from the defining recurrence for the sequence.
Let $\a_n$ denote this set of paths. 
\begin{center}
\begin{pspicture}(-5,-1.6)(15,2)
\psset{unit=.4cm}
\rput(-0.5,1.6){\footnotesize \mag{0}}
\rput(2.3,1.5){\footnotesize \mag{2}}
\rput(3.7,1.5){\footnotesize \mag{4}}
\rput(8.3,1.5){\footnotesize \mag{8}}
\rput(16.3,1.5){\footnotesize \mag{16}}
\rput(-0.9,3){\footnotesize \mag{first return\:\raisebox{-2mm}{$\searrow$}}}
\rput(12,1.3){\mag{\footnotesize $\mbf{\uparrow}$}}
\rput(12.6,0.4){\mag{\footnotesize $x$-axis, or \gl}} 
\psdots(0,2)(1,1)(2,2)(3,3)(4,2)(5,1)(6,0)(7,1)(8,2)(9,3)(10,4)(11,5)(12,4)(14,4)(15,3)(16,2)
\psline(0,2)(1,1)(2,2)(3,3)(4,2)(5,1)(6,0)(7,1)(8,2)(9,3)(10,4)(11,5)(12,4)(14,4)(15,3)(16,2)
\psline[linecolor=magenta,linestyle=dotted](0,2)(16,2)

\rput(8,-1.7){\small A path in $\a_8$} 
\rput(8,-3.9){Figure 1} 

\end{pspicture}
\end{center}

It is not too hard to find the \gf for  $\a_n$ using the ``first return'' decomposition and a little recursion involving both nonnegative paths and the analogous \gfs when ``\,$y=2$\,'' in the defining condition is replaced by ``\,$y=1$\,'' and by ``\,$y=0$\,'', and then observe that the \gf for $\a_n$ coincides with that for A111279. This answers Plewe's question in the affirmative. But it's nicer to give a bijective proof.

By definition, the $(n+1)$-th term of A111279 counts permutations of $[n+1]$ avoiding the 
three patterns \{3241, 3421, 4321\}. 
These permutations are in bijection \cite{weaksort} with the set $\b_n$ of Schr\"{o}der paths of size $n$ \emph{with at most one peak per component}. 
Recall that a \emph{Schr\"{o}der path} is 
a nonnegative Grand Schr\"{o}der path, that is, one that never dips below the $x$-axis, and the interior vertices on the $x$-axis split a nonempty path that ends on the $x$-axis into its components.
\begin{center}
\begin{pspicture}(-5,-1.5)(15,1.2)
\psset{unit=.4cm}
\psdots(0,0)(2,0)(3,1)(4,0)(5,1)(7,1)(8,0)(9,1)(10,2)(12,2)(13,3)(14,2)(15,1)(16,0)
\psline(0,0)(2,0)(3,1)(4,0)(5,1)(7,1)(8,0)(9,1)(10,2)(12,2)(13,3)(14,2)(15,1)(16,0)
\psline[linecolor=green,linestyle=dotted](0,0)(16,0)

\rput(8,-1.6){\small A path in $\b_8$ with 4 components} 
\rput(8,-3.8){Figure 2} 

\end{pspicture}
\end{center}

\section{The bijection} 

\vspace*{-2mm}

We give a simple bijection of the cut-and-paste type from $\a_n$ to $\b_n$ that preserves 
components. Thus it suffices to define our mapping on indecomposable (1-component) paths 
in $\a_n$. We refer to the horizontal line joining the terminal points of a path in $\a_n$ 
or $\b_n$ as \emph{\gl}, GL for short, to eliminate the need for coordinate axes.
Since a path in $\a_n$ contains no flatsteps at \gl, an indecomposable path in $\a_n$ 
cannot consist of a single flatstep and so lies entirely above or entirely below \gl. 

If entirely below (and hence contains no flatsteps at all), flip it 
over the $x$-axis and replace all $UD$s (peaks) by $F$s to get an indecomposable 
Schr\"{o}der $n$-path with \emph{no} peak (Figure 3). 
\begin{center}
\begin{pspicture}(-0.5,-3.5)(35,1.5)
\psset{unit=.4cm}
\psdots(0,0)(1,-1)(2,-2)(3,-1)(4,-2)(5,-3)(6,-2)(7,-1)(8,0)
\psdots(15,0)(16,1)(17,2)(18,1)(19,2)(20,3)(21,2)(22,1)(23,0)
\psdots(30,0)(31,1)(33,1)(34,2)(36,2)(37,1)(38,0)
\psline(0,0)(1,-1)(2,-2)(3,-1)(4,-2)(5,-3)(6,-2)(7,-1)(8,0) 
\psline(15,0)(16,1)(17,2)(18,1)(19,2)(20,3)(21,2)(22,1)(23,0)
\psline(30,0)(31,1)(33,1)(34,2)(36,2)(37,1)(38,0)
\psline[linecolor=green,linestyle=dotted](0,0)(8,0)
\psline[linecolor=green,linestyle=dotted](15,0)(23,0)
\psline[linecolor=green,linestyle=dotted](30,0)(38,0)

\rput(11.5,1.2){\footnotesize flip} 
\rput(11.5,0){$\longrightarrow$}

\rput(26.5,2.4){\footnotesize flatten}
\rput(26.5,1.2){\footnotesize all peaks} 
\rput(26.5,0){$\longrightarrow$}

\rput(4,-4.8){\small path in $\a_4$} 
\rput(4,-6){\small indecomposable, below GL}

\rput(34,-2.3){\small path in $\b_4$} 
\rput(34,-3.5){\small indecomposable, no peaks}
\rput(19.5,-8.5){Figure 3}
\end{pspicture}
\end{center}

If entirely above, follow the sequence of operations illustrated in Figure 4 below to get an indecomposable Schr\"{o}der $n$-path with \emph{exactly one} peak.

\begin{center}
\begin{pspicture}(0,-2)(15,2)
\psset{unit=.4cm}
\psdots(-1,-1)(0,0)(1,1)(2,2)(3,1)(4,0)(5,1)(7,1)(8,2)(9,3)(10,2)(11,3)(12,2)(13,1)(14,2)(15,1)(16,0)(17,1)(18,0)(19,1)(20,2)(21,1)(22,0)(23,-1)
\psline(-1,-1)(0,0)(2,2)(4,0)(5,1)(7,1)(9,3)(10,2)(11,3)(13,1)(14,2)(16,0)(17,1)(18,0)(20,2)(22,0)(23,-1)
\psline[linecolor=green,linestyle=dotted](-1,-1)(23,-1)
\rput(32,0){$\longrightarrow$}
\rput(32,2.3){\textrm{{\footnotesize delete first and}}}
\rput(32,1.2){\textrm{{\footnotesize last steps}}}
\rput(11,-2.3){\small indecomposable path in $\a_{12}$}
\rput(11,-3.5){\small above GL, all $F$s at height 2} 
\end{pspicture}
\end{center}
\begin{center}
\begin{pspicture}(0,0)(15,2)
\psset{unit=.4cm}
\psdots(0,0)(1,1)(2,2)(3,1)(4,0)(5,1)(7,1)(8,2)(9,3)(10,2)(11,3)(12,2)(13,1)(14,2)(15,1)(16,0)(17,1)(18,0)(19,1)(20,2)(21,1)(22,0)
\psline(0,0)(2,2)(4,0)(5,1)(7,1)(9,3)(10,2)(11,3)(13,1)(14,2)(16,0)(17,1)(18,0)(20,2)(22,0)
\psline[linecolor=green,linestyle=dotted](0,0)(22,0)

\rput(32,0){$\longrightarrow$}
\rput(32,2.1){\textrm{{\footnotesize replace each $F$ with a $DU$,}}}
\rput(32,1.1){\textrm{{\footnotesize hiliting new vertices at GL}}}
\end{pspicture}
\end{center}
\begin{center}
\begin{pspicture}(0,0)(15,2)
\psset{unit=.4cm}
\psdots(0,0)(1,1)(2,2)(3,1)(4,0)(5,1)(6,0)(7,1)(8,2)(9,3)(10,2)(11,3)(12,2)(13,1)(14,2)(15,1)(16,0)(17,1)(18,0)(19,1)(20,2)(21,1)(22,0)
\psline(0,0)(2,2)(4,0)(5,1)(6,0)(7,1)(9,3)(10,2)(11,3)(13,1)(14,2)(16,0)(17,1)(18,0)(20,2)(22,0)
\psline[linecolor=green,linestyle=dotted](0,0)(22,0)
 \psdots[linewidth=1mm](6,0)
\rput(32,0){$\longrightarrow$}
\rput(32,2.3){\textrm{{\footnotesize flip first component as well as those}}}
\rput(32,1.2){\textrm{{\footnotesize (if any) that start at a hilited vertex}}}
\end{pspicture}
\end{center}
\begin{center}
\begin{pspicture}(0,-2)(15,2)
\psset{unit=.4cm}
\psdots(0,0)(1,-1)(2,-2)(3,-1)(4,0)(5,1)(6,0)(7,-1)(8,-2)(9,-3)(10,-2)(11,-3)(12,-2)(13,-1)(14,-2)(15,-1)(16,0)(17,1)(18,0)(19,1)(20,2)(21,1)(22,0)
\psline(0,0)(1,-1)(2,-2)(3,-1)(4,0)(5,1)(6,0)(7,-1)(8,-2)(9,-3)(10,-2)(11,-3)(12,-2)(13,-1)(14,-2)(15,-1)(16,0)(17,1)(18,0)(19,1)(20,2)(21,1)(22,0)
\psline[linecolor=green,linestyle=dotted](0,0)(22,0)
\rput(32,-2.1){$\longrightarrow$}
\rput(32,1.3){\textrm{{\footnotesize now $V_1$ is leftmost of the lowest vertices,}}}
\rput(32,0.2){\textrm{{\footnotesize and $V_2$ terminates the last $U$ ending at GL;}}}
\rput(32,-0.9){\textrm{{\footnotesize interchange segments $O$ to $V_1$ and $V_1$ to $V_2$}}}

\rput(-0.2,.6){\textrm{\footnotesize $O$}}
\rput(9,-3.8){\textrm{\footnotesize $V_1$}}
\rput(16.3,-.7){\textrm{\footnotesize $V_2$}}
\end{pspicture}
\end{center}
\begin{center}
\begin{pspicture}(0,-.5)(15,2)
\psset{unit=.4cm}
\psline(0,0)(1,1)(2,0)(3,1)(4,2)(5,1)(6,2)(7,3)(8,2)(9,1)(10,2)(11,3)(12,4)(13,3)(14,2)(15,1)(16,0)(17,1)(18,0)(19,1)(20,2)(21,1)(22,0)

\psdots(0,0)(1,1)(2,0)(3,1)(4,2)(5,1)(6,2)(7,3)
\psdots(8,2)(9,1)(10,2)(11,3)(12,4)(13,3)(14,2)(15,1)(16,0)(17,1)(18,0)(19,1)(20,2)(21,1)(22,0)
\psline[linecolor=green,linestyle=dotted](0,0)(22,0)

\rput(7,3.8){\textrm{\footnotesize $V_2$}}

\rput(32,0){$\longrightarrow$}
\rput(32,2.2){\textrm{{\footnotesize flatten all peaks }}} 
\rput(32,1.1){\textrm{{\footnotesize except the one at $V_2$}}}
\end{pspicture}
\end{center}
\begin{center}
\begin{pspicture}(0,-.5)(15,2)
\psset{unit=.4cm}
\psline[linecolor=green,linestyle=dotted](0,0)(22,0)
\psline(0,0)(2,0)(3,1)(5,1)(6,2)(7,3)(8,2)(9,1)(10,2)(11,3)(13,3)(14,2)(15,1)(16,0)(18,0)(20,2)(22,0)

\psdots(0,0)(2,0)(3,1)(5,1)(6,2)(7,3)
\psdots(8,2)(9,1)(10,2)(11,3)(13,3)(14,2)(15,1)(16,0)(18,0)(19,1)(20,2)(21,1)(22,0)

\rput(32,0){$\longrightarrow$}
\rput(32,1.1){\textrm{{\footnotesize prepend $U$, append $D$}}}

\end{pspicture}
\end{center}
\begin{center}
\begin{pspicture}(0,-1.5)(15,2)
\psset{unit=.4cm}
\psline[linecolor=green,linestyle=dotted](-1,-1)(23,-1)
\psline(-1,-1)(0,0)(2,0)(3,1)(5,1)(6,2)(7,3)(8,2)(9,1)(10,2)(11,3)(13,3)(14,2)(15,1)(16,0)(18,0)(20,2)(22,0)(23,-1)

\psdots(-1,-1)(0,0)(2,0)(3,1)(5,1)(6,2)(7,3)
\psdots(8,2)(9,1)(10,2)(11,3)(13,3)(14,2)(15,1)(16,0)(18,0)(19,1)(20,2)(21,1)(22,0)(23,-1)

\rput(25,1){$\leftarrow$}
\rput(33,1.6){\textrm{{\small indecomposable Schr\"{o}der path}}}
\rput(33,.4){\textrm{{\small in $\b_{12}$ with exactly one peak }}}
\rput(20,-3){Figure 4}

\end{pspicture}
\end{center}

All the steps are reversible, and we have the desired bijection.

Exercise. Check that the $\a_8$-path in Figure 1 corresponds to the $\b_8$-path in Figure 2 under this bijection.

\noindent \htmladdnormallink{Department of Statistics}{http://www.stat.wisc.edu/}, University of Wisconsin-Madison  \\
\noindent 1300 University Ave, Madison, WI \ 53706-1532  \\
{\bf callan@stat.wisc.edu}  \\

\end{document}